\documentclass[11pt]{article}\include{ddefs}
\usepackage{graphics,pifont,array}
\usepackage{graphicx}
\usepackage{endnotes}
\begin{document}

\begin{titlepage}

\begin{center}

{\LARGE{        \mbox{   }                          \\
                \mbox{   }                          \\
                \mbox{   }                          \\
                \mbox{   }                          \\
                \mbox{   }                          \\
                \mbox{   }                          \\
                \mbox{   }                          \\
                \mbox{   }                          \\
  {\bf{Homogeneous and Inhomogeneous Integral Formulations of \\
       Nonrelativistic Potential Scattering}}  \\
               \mbox{    }                          \\
               \mbox{    }                          \\
              J. A. Grzesik                         \\
           Allwave Corporation                      \\
        3860 Del Amo Boulevard                      \\
                Suite 404                           \\
           Torrance, CA 90503                       \\
                \mbox{    }                         \\
        +1 (818) 749-3602                           \\ 
            jan.grzesik@hotmail.com                 \\
              \mbox{     }                          \\
              \mbox{     }                          \\  
             June 05, 2020                              }   }

\end{center}

\end{titlepage}

\pagenumbering{roman}

\setcounter{page}{2}

\vspace*{+2.725in}

\begin{abstract}

	\parindent=0.245in
	
	Advantage is taken of the arbitrariness in energy reference to consider anew
	integral transcriptions of Schr{\"{o}}dinger's equation in the presence of
	potentials which at infinity acquire constant, nonvanishing values.  It is found
	possible to present for the probability amplitude $\psi$ a linear integral equation
	which is entirely devoid of explicit reference to the wave function incident from
	infinity, and thus differs markedly from the prevailing inhomogeneous formulation.
	Identity of the homogeneous equation with an inhomogeneous statement which is at
	the same time available is affirmed in general terms with the aid of the
	Fourier transformation, and is then still further reinforced by application of both
	formalisms to the particular example of a spherical potential barrier/well.
	Identical, closed-form outcomes are gotten in each case for wave function eigenmode
	expansion coefficients on both scatterer interior and exterior.  Admittedly, the
	solution procedure is far simpler in the inhomogeneous setting, wherein it exhibits
	the aspect of a direct, leapfrog advance, unburdened by any implicit algebraic entanglement.
    By contrast, the homogeneous path, of considerably greater length, insists, at each
    mode index,	upon an exterior/interior coefficient entanglement, an entanglement which,
    happily, is no more severe than	that of a non-singular two-by-two linear system.  Each
    such two-by-two linear system reproduces of course the output already gotten under
    the inhomogeneous route, and is indeed identical to the two-by-two system encountered
    during the routine procedure wherein continuity is demanded at the barrier/well interface
    of both $\psi$ and its radial derivative.

\end{abstract}

\newpage

\parindent=0.5in

\pagenumbering{arabic}

\setcounter{page}{1}

\pagestyle{myheadings}

\mbox{   }

\markright{J. A. Grzesik \\ homogeneous and inhomogeneous formulations of nonrelativistic potential scattering}

\section{Introduction}
     A transcription of the time-independent Schr{\"{o}}dinger differential equation in the form of a
linear inhomogeneous integral equation represents the cornerstone of most discussions of quantum
mechanical potential scattering.  In these discussions the reference of energy is implicitly understood
to vanish at infinity, $V_{\infty}=0.$  The potentials $V(\mbox{\boldmath$r$})$ are assumed to be of bounded
support or else to approach their asymptotic null value with sufficient rapidity to assure convergence
of all integrals involved.  Suppose however that the reference of energies is displaced so that
$V_{\infty}\neq 0.$  We wish to indicate that under these circumstances one has available a choice
of integral transcriptions, one homogeneous while the other inhomogeneous, which are entirely
equivalent to each other.  The inhomogeneous equation follows from a comparison of all energies
with $V_{\infty}.$  It requires that the incident particle wave be explicitly displayed and permits
an easy transition to the limit in which $V_{\infty}\rightarrow 0.$  Construction of the homogeneous
equation, on the other hand, relies critically upon the fact that $V_{\infty}\neq 0,$ but otherwise
demands no explicit provision for the incident field.  Equivalence of these formulations is difficult
to establish directly in configuration space, but it does yield to an easy demonstration when the integral
equations are first subjected to a Fourier transformation.  That equivalence is then further
illustrated by an application of the integral equations to the simple concrete example involving
scattering by a spherically symmetric potential barrier/well.

\section{A homogeneous/inhomgeneous (h/i) integral equation duality}

   We represent the potential $V(\mbox{\boldmath$r$})$ in the form
\begin{equation}
V(\mbox{\boldmath$r$}) = v(\mbox{\boldmath$r$})+V_{\infty}\,,
\end{equation}
with $V_{\infty}\neq 0$ and $v(\mbox{\boldmath$r$})\rightarrow 0$ as $\mbox{\boldmath$r$}$ grows
without bound in all directions.  A plane wave energy-momentum eigenstate
\begin{equation}
\psi_{inc}(\mbox{\boldmath$r$}) =  \exp(i\mbox{\boldmath$k\cdot r$})\,,
\end{equation}
with
\begin{equation}
k^{2}  =  \frac{2m}{\hbar^{2}}\left(\rule{0mm}{4mm}E-V_{\infty}\right) \, > \, 0\,,
\end{equation}
is incident from infinity upon the scattering potential $v(\mbox{\boldmath$r$}).$  Symbols $m,$ $E,$
and $\mbox{\boldmath$k$}$ have their usual meanings as the particle mass, total energy, and vector
momentum, the latter being measured in units of Planck's constant $\hbar.$

      In order to place the ensuing homogeneous formulation in a proper perspective, we recapitulate
briefly a statement of the inhomogeneous integral equation.  This latter depends upon explicit
decomposition of the wave function $\psi$ into incident and scattered contributions
\begin{equation}
\psi  =  \psi_{inc}+\psi_{scatt}\,.
\end{equation}
One obtains for $\psi_{scatt}$ the equation
\newpage
\mbox{    }
\mbox{    }
\newline
\newline
\begin{equation}
\left(\rule{0mm}{4mm}\nabla^{2}+k^{2}\right)\psi_{scatt}  =  \frac{2m}{\hbar^{2}}v\psi\,,
\end{equation}               
with the transcription
\begin{equation}
\psi_{scatt}(\mbox{\boldmath$r$}) = \frac{2m}{\hbar^{2}}\int G(|\mbox{\boldmath$r-r'$}|)
             v(\mbox{\boldmath$r$})\psi(\mbox{\boldmath$r'$})d\mbox{\boldmath$r'$}\,,
\end{equation}
wherein
\begin{equation}
G(r)  =  -\frac{e^{ikr}}{4\pi r}
\end{equation}
is the Green's function for the Helmholtz equation which involves expanding waves alone.  The composite
wave function $\psi$ from (4) therefore satisfies
\begin{equation}
\psi(\mbox{\boldmath$r$}) = \psi_{inc}(\mbox{\boldmath$r$})+\frac{2m}{\hbar^{2}}\int G(|\mbox{\boldmath$r-r'$}|)
      v(\mbox{\boldmath$r'$})\psi(\mbox{\boldmath$r'$})d\mbox{\boldmath$r'$}\,,      
\end{equation}
the integrals in (8) and elsewhere being taken over all of space and differential $d\mbox{\boldmath$r$}$
being a standard shorthand for the product $dxdydz.$  To this point the sole modification of the usual integral
formulation reposes in the representation (3) which follows from the displacement $V_{\infty}$ in energy
reference.  Transition to the limit $V_{\infty}\rightarrow 0$ in (8) is automatic, provided of course that energy
eigenvalue $E$ shifts accordingly so as to keep the difference on the right in (3) fixed, an invariance which
respects the fact that the operational datum here consists entirely of the incoming kinetic energy and direction
of particle flight.

     Erection of a homogeneous integral equation for $\psi$ becomes possible if at the outset one forgoes the
decomposition (4).  The Schr{\"{o}}dinger equation in the form
\begin{equation}
\left(\rule{0mm}{4mm}\nabla^{2}+k_{0}^{2}\right)\psi = \frac{2m}{\hbar^{2}}V\psi\,,
\end{equation}
with
\begin{equation}
k_{0}^{2}  =  \frac{2m}{\hbar^{2}}E\,,
\end{equation}
entails a source term which, unlike that in (5), persists to infinity.  From (9) one therefore obtains
\begin{equation}
\psi(\mbox{\boldmath$r$})  =  \frac{2m}{\hbar^{2}}\int G_{0}(|\mbox{\boldmath$r-r'$}|)
V(\mbox{\boldmath$r'$})\psi(\mbox{\boldmath$r'$})d\mbox{\boldmath$r'$}\,,
\end{equation}
as the homogeneous counterpart of (8).  In (11) the suffix zero serves to indicate that $G_{0}$ is obtained from
(7) by the replacement $k\rightarrow k_{0}.$\footnote{We hasten to remark that, on physical grounds, a positivity requirement
constrains only (3) on its right, and that for a sufficiently negative choice of $V_{\infty}$ there is nothing to
prevent having total energy $E$ itself negative.  Should that occur, we become required of course to take $k_{0}$
as positive imaginary, {\em{viz.,}} $k_{0}=i\sqrt{2m|E|/\hbar^{2}\,},$ so as to retain integral convergence at
infinity.}
\section{Demonstration of equivalence}
All reasonable attempts to supply a direct, configuration space demonstration of identity for (8) and (11)
encounter what seem to be insuperable difficulties.\mbox{ }  The required identity, however, emerges quite simply
\newpage
\mbox{   }
\newline
\newline
\newline
under Fourier transformation, which for definiteness is taken here in the asymmetric form
\begin{equation}
F(\cdots)  =  \frac{1}{(2\pi)^{3}}\int \exp(-i\mbox{\boldmath$K\cdot r$}) \cdots\, d\mbox{\boldmath$r$}\,,
\end{equation}
its outcome being designated by a tilde placed above and transform argument $\mbox{\boldmath$K$}$ replacing its
heritage antecedent $\mbox{\boldmath$r$}.$  The underlying unknown then becomes the scattered wave function
transform ${\tilde{\psi}}(\mbox{\boldmath$K$}).$  Identity of the transformed integral equations for
${\tilde{\psi}}(\mbox{\boldmath$K$})$ is understood to constitute a demonstration of equivalence for (8) and
(11).

     Fourier transformation is immediate if one exploits the transformation law for spatial convolution
products and takes into account the subordinate results

\vspace{-6mm}
{\large{
\begin{equation}
\left.\rule{0mm}{1.1cm}\begin{array}{rcl}
F\left(\rule{0mm}{3.5mm}(r\exp(-ikr))^{-1}\right) & = &
     \left\{\rule{0mm}{3.5mm} 2\pi^{2}(K^{2}-k^{2})\right\}^{-1} \\   
                      &   &                                   \\
F\left(\rule{0mm}{3.5mm}\psi_{inc}(\mbox{\boldmath$r$})\right) & = & \delta(\mbox{\boldmath$K-k$})
       \end{array} \right\}\,,                     
\end{equation} } }
\parindent=0in

\vspace{-6mm}   
with $\delta(\mbox{\boldmath$K-k$})$ being the Dirac delta function in three dimensions.
In the neighborhood of the shell $K=k$ the first of the structures in (13) assumes the usual
form\footnote{The plus sign prefacing Dirac's delta is associated with our implicit understanding that
wave number $k$ is to be regarded as the limit of $k+i\epsilon$ when $\epsilon\downarrow0+,$
a standard, formal device contrived to assure integral convergence at infinity.  One then finds
that $(k+i\epsilon)^{2}=k^{2}-\epsilon^{2}+2i\epsilon k$ likewise falls above the real axis,
circumstance which ultimately leads to a plus sign being attached to Dirac's delta in (14).}
\parindent=0.5in  	
\begin{equation}
\frac{1}{2\pi^{2}(K^{2}-k^{2})} =  \frac{1}{2\pi^{2}}\left(\rule{0mm}{6mm}\frac{P}{K^{2}-k^{2}}+
                        \frac{i\pi}{2k}\delta(K-k)\right)\,
\end{equation}
wherein $P$ represents Cauchy's principal value whereas $\delta(K-k)$ is the one-dimensional analogue of
$\delta(\mbox{\boldmath$K-k$}).$  Abbreviate also by writing
\begin{eqnarray}
S(\mbox{\boldmath$K$}) & = & F(v\psi_{scatt}) \nonumber \\
     & = & \int {\tilde{v}}(\mbox{\boldmath$K-K'$}){\tilde{\psi}}_{scatt}(\mbox{\boldmath$K'$})d\mbox{\boldmath$K'$}\,.
\end{eqnarray}
Then on the basis of (8) one obtains
\begin{equation}
{\tilde{\psi}}_{scatt}(\mbox{\boldmath$K$}) = - \frac{2m}{\hbar^{2}(K^{2}-k^{2})}\left(\rule{0mm}{6mm}
       {\tilde{v}}(\mbox{\boldmath$K-k$})+S(\mbox{\boldmath$K$})\right)\,,
\end{equation}
an inhomogeneous integral equation in its own right, but now in reciprocal, $\mbox{\boldmath$K$}$ space.

     Analogous steps undertaken in the context of (11) then
yield\footnote{The divisions by $K^{2}-k_{0}^{2}$ appearing on the right in (17) and on the left
in (18) are sufficiently tame to cause no real trouble:  if $k_{0}^{2}<0,$ there is clearly no
problem, whereas if $k_{0}^{2}>0,$ then a vanishing, positively signed imaginary accompaniment
to $k_{0}$ is understood as before with $k,$ leading to an obvious analogue of (14).}
\begin{equation}
\delta(\mbox{\boldmath$K-k$})+{\tilde{\psi}}_{scatt}(\mbox{\boldmath$K$})  = 
-\frac{2m}{\hbar^{2}(K^{2}-k_{0}^{2})}\left(\rule{0mm}{6mm}
{\tilde{v}}(\mbox{\boldmath$K-k$})+S(\mbox{\boldmath$K$})+V_{\infty}
\left\{\rule{0mm}{5mm}\delta(\mbox{\boldmath$K-k$})+{\tilde{\psi}}_{scatt}(\mbox{\boldmath$K$})\right\}\right)\,.   
\end{equation}
\newline
\newline
\newline
Here
\begin{eqnarray}
\frac{2mV_{\infty}}{\hbar^{2}(K^{2}-k_{0}^{2})}\,\delta(\mbox{\boldmath$K-k$}) & = &
          \frac{2mV_{\infty}}{\hbar^{2}(k^{2}-k_{0}^{2})}\,\delta(\mbox{\boldmath$K-k$}) \nonumber \\
     & = &  -\delta(\mbox{\boldmath$K-k$})    
\end{eqnarray}
when note is taken of (3) and (10), so that the delta functions in (17) cancel identically.  A
simple transposition of terms suffices next to display (17) in the form of (16).  This concludes the
formal demonstration of equivalence.

    One small trace of this equivalence is already afforded by noting that $\psi_{inc}$ as given in (2) 
itself satisfies (11) when that latter is specialized by setting $v=0.$  The first of the Fourier
transforms (13) is an essential ingredient in this sudsidiary aspect of equivalence.
\section{The limit $V_{\infty}\rightarrow 0,$ $k_{0}\rightarrow k$}
    The equivalence of (8) and (11) just now established provides an {\em{a priori}} assurance that the limit
of (11) as $V_{\infty}\rightarrow 0,$ $k_{0}\rightarrow k$ must coalesce with that of (8).  Inasmuch as in this
limit there necessarily occurs a transition from homogeneous to inhomogeneous integral equations, a
more direct discussion is of some independent interest.  It bears repeating perhaps that energy
eigenvalue $E$ slides up or down in unison with $V_{\infty},$ subject only to the overriding
demand that the right hand side of (3) remain invariant.

     Consider (11) written out in full as
\begin{eqnarray}
\psi(\mbox{\boldmath$r$}) & = & \frac{2mV_{\infty}}{\hbar^{2}}
       \int G_{0}(|\mbox{\boldmath$r-r'$}|)\psi_{inc}(\mbox{\boldmath$r'$})d\mbox{\boldmath$r'$} \nonumber \\
     &    & \rule{1cm}{0mm}+ \frac{2mV_{\infty}}{\hbar^{2}}
       \int G_{0}(|\mbox{\boldmath$r-r'$}|)\psi_{scatt}(\mbox{\boldmath$r'$})d\mbox{\boldmath$r'$}  \\
     &    & \rule{2cm}{0mm} + \frac{2m}{\hbar^{2}}
     \int G_{0}(|\mbox{\boldmath$r-r'$}|)v(\mbox{\boldmath$r'$})\psi(\mbox{\boldmath$r'$})d\mbox{\boldmath$r'$}\,.\nonumber  
\end{eqnarray}
From the concluding remark of the previous section the first term on the right is simply $\psi_{inc}(\mbox{\boldmath$r$}),$
independently of $V_{\infty}.$  All that remains to display identity of limits for (8) and (11) therefore is to show
that the second line vanishes together with $V_{\infty}.$  This behavior in turn is guaranteed as soon as freedom from
divergence is exhibited for the integral
\begin{equation}
\int G_{0}(|\mbox{\boldmath$r-r'$}|)\psi_{scatt}(\mbox{\boldmath$r'$})d\mbox{\boldmath$r'$}\,.
\end{equation}

   Now the only possible source of divergence for (20) is integration in the asymptotic region
$r'\rightarrow\infty.$  Over this domain one writes as always,
with $\mbox{\boldmath${\hat{r}}'$}=\mbox{\boldmath$r'$}/r',$
\begin{equation}
G_{0}(|\mbox{\boldmath$r-r'$}|) \approx  -\frac{e^{ik_{0}r'}}{4\pi r'}
     \, \exp\left(\rule{0mm}{3.25mm}\!-ik_{0}\mbox{\boldmath${\hat{r}}'\cdot r$}\right)\,,
\end{equation}
\newpage
\mbox{  }
\newline
\newline
\newline
while for $\psi_{scatt}$ the corresponding representation
\begin{equation}
\psi_{scatt}(\mbox{\boldmath$r'$}) \rule{3mm}{0mm} \approx_{\atop{\rule{-5.5mm}{3mm}r'\rightarrow\infty}} 
 -\left(\frac{m}{2\pi\hbar^{2}}\right)\!\left(\frac{e^{ikr'}}{r'}\right)
 \int\exp\left(\rule{0mm}{3.25mm}\!-ik\mbox{\boldmath${\hat{r}}'\cdot r''$}\right)v(\mbox{\boldmath$r''$})
      \psi(\mbox{\boldmath$r''$})d\mbox{\boldmath$r''$}
\end{equation}
provided by (8) is employed.  In view of the established equivalence of (8) and (11) this appeal to (8) is entirely
legitimate.  The convergence of (20) is thus seen to depend upon that of
\begin{equation}
\int_{0}^{\,\infty} e^{i(k_{0}+k)r'}dr' =  \frac{i}{k_{0}+k} 
\end{equation}
and hence is amply assured by virtue of all previous remarks concerning value ranges admissible to
$k_{0}$ and $k.$  Accordingly, (20) is indeed a convergent integral and the limits of (8) and (11) are
indeed identical.
\section{A simple example:  spherical potential barrier/well}
    While the pre{\"{e}}minent utility of integral scattering formulations rests upon the easy access which they
provide to iterative approximation schemes, it is possible also to base upon them exact solutions of simple
scattering problems.  In such approaches one avoids all questions of continuity at geometric boundaries
separating regions of differing potential, and requires merely that the integral equations be satisfied both
interior and exterior to such boundaries.  This requirement is completely adequate to determine all expansion
coefficients when mode decompositions are employed.

    Consider for example a spherical region $0\leq r\leq a$ throughout which $V$ equals a constant
$V_{1}\neq V_{\infty}.$  For $r>a$ the potential assumes the uniform value $V_{\infty}\neq 0.$  Otherwise
put, $v(\mbox{\boldmath${r}$})=V_{1}-V_{\infty}$ whenever $0\leq r \leq a$ and is zero otherwise.
Positive $v$ betokens a repulsive barrier, negative $v$ an attractive well. 
Then for $r>a$ one writes in terms of Legendre polynomials $P_{l}$ and spherical Hankel functions
$h^{(1)}_{l}$ in standard notation
\begin{equation}
\psi_{scatt}(\mbox{\boldmath$r$})  = \sum_{l=0}^{\infty}A^{>}_{l}h^{(1)}_{l}(kr)
    P_{l}(\mbox{\boldmath${\hat{k}}\cdot{\hat{r}}$})\,,
\end{equation}
while for $0\leq r \leq a$ there holds for the total wave field a similar decomposition
\begin{equation}
\psi(\mbox{\boldmath$r$}) =  \sum_{l=0}^{\infty}A^{<}_{l}j_{l}(k_{1}r)
P_{l}(\mbox{\boldmath${\hat{k}}\cdot{\hat{r}}$})
\end{equation}
in terms of spherical Bessel functions $j_{l}$ governed by an interior propagation constant
\begin{equation}
k_{1}  =  \sqrt{\frac{2m}{\hbar^{2}}\left(\rule{0mm}{4mm}E-V_{1}\right)\,}\,.
\end{equation}
Classical inaccessibility of the region $0\leq r \leq a$ corresponds as always to a purely imaginary value for
$k_{1}.$  Our task is to determine the interior/exterior expansion coefficients $\{A^{>}_{l}\}_{l=0}^{\infty}$
and $\{A^{<}_{l}\}_{l=0}^{\infty}$ on the basis of Eqs. (8) and (11).
\newpage
\mbox{   }
\newline

   One employs also the standard representations
\begin{equation}
\psi_{inc}(\mbox{\boldmath${r}$})  =   \sum_{l=0}^{\infty}i^{l}(2l+1)j_{l}(kr)
           P_{l}(\mbox{\boldmath${\hat{k}}\cdot{\hat{r}}$})\,,
\end{equation}
\begin{equation}
G(|\mbox{\boldmath$r-r'$}|)  =  \frac{k}{4\pi i}\sum_{l=0}^{\infty}(2l+1)j_{l}(kr_{<})h^{(1)}_{l}(kr_{>})
       P_{l}(\mbox{\boldmath${\hat{r}}\cdot{\hat{r}}'$})\,,
\end{equation}   
and similarly for $G_{0}.$  We adhere to common practice by writing $r_{>}$ and $r_{<}$ respectively for
the greater and lesser of $r$ and $r',$ which is to say $r_{>}=\max(r,r')=(r+r'+|r-r'|)/2,$
$r_{<}=\min(r,r')=(r+r'-|r-r'|)/2.$  In (27), as opposed to what is both physically plausible and mathematically
desirable in (24) and (28), wave number $k$ must be taken strictly real so as to assure an incoming wave field that remains
bounded at infinity, both fore and aft.
\vspace{-2mm}  
\subsection{Solution based on the inhomogeneous form (8)}
\vspace{-2mm}
    We examine first an application of the inhomogeneous form (8) to a determination of the expansion coefficients
$\{A^{>}_{l}\}_{l=0}^{\infty}$ and $\{A^{<}_{l}\}_{l=0}^{\infty}.$  The requirement that the integral equation hold
true for $r>a$ yields
\begin{eqnarray}
\lefteqn{
\sum_{l=0}^{\infty}A^{>}_{l}h_{l}^{(1)}(kr)P_{l}(\mbox{\boldmath${\hat{k}}\cdot{\hat{r}}$})  = } \nonumber  \\
    &    &  \frac{mk}{2\hbar^{2}\pi i}\left(\rule{0mm}{4mm}V_{1}-V_{\infty}\right)\sum_{l,n=0}^{\infty}A^{<}_{l}(2n+1)
 	h^{(1)}_{n}(kr)\int_{r'\leq a}P_{l}(\mbox{\boldmath${\hat{k}}\cdot{\hat{r}}'$})
 	P_{n}(\mbox{\boldmath${\hat{r}}\cdot{\hat{r}}'$})j_{l}(k_{1}r')j_{n}(kr')d\mbox{\boldmath$r'$}\,. 
\end{eqnarray}
Integrations over the angular c{\"{o}}ordinates $(\vartheta',\varphi')$ of $\mbox{\boldmath${r'}$}$ can be effected by appealing to the addition theorem for Legendre polynomials.  It is then
easily found that the matrix of integrals in (29) is diagonal in its indices $l$ and $n,$ with the element at spot
$(l,l)$ proportional to $P_{l}(\mbox{\boldmath${\hat{k}}\cdot{\hat{r}}$}).$  Both members of (29), left and right, assume
the form of single series in the functions $h_{l}^{(1)}(kr)P_{l}(\mbox{\boldmath${\hat{k}}\cdot{\hat{r}}$})$ and so
force equality of corresponding coefficients.  Taking note also of a standard (Lommel) quadrature formula for Bessel function
pairs having a common index, and of the relation
\begin{equation}
V_{1}-V_{\infty} = \frac{\hbar^{2}}{2m}\left(\rule{0mm}{3.00mm}k^{2}-k_{1}^{2}\right)\,,
\end{equation}
\begin{equation}
A^{>}_{l} = \frac{\pi a}{2i}\sqrt{\frac{k}{k_{1}}\,}\left(\rule{0mm}{6mm}
                  k_{1}J_{l+\frac{1}{2}}'(k_{1}a)J_{l+\frac{1}{2}}(ka)
                       -kJ_{l+\frac{1}{2}}(k_{1}a)J_{l+\frac{1}{2}}'(ka)\right)A^{<}_{l}\,,
\end{equation}
with spherical Bessel functions restored to their half-index antecedents.  The spherical weight factor $r'^{\,2}$
in the quadratures on the right in (29) has of course been taken into account. 

     To enforce the validity of (8) for $0\leq r \leq a$ one writes
\begin{eqnarray}
\lefteqn{
\sum_{l=0}^{\infty}\left(\rule{0mm}{5mm}A^{<}_{l}j_{l}(k_{1}r)-i^{l}(2l+1)j_{l}(kr)\right)
                   P_{l}(\mbox{\boldmath${\hat{k}}\cdot{\hat{r}}$})  = }  \nonumber   \\ 
    &     &  k\left(\frac{k^{2}-k_{1}^{2}}{4\pi i}\right)\sum_{l,n=0}^{\infty}A^{<}_{l}(2n+1) 
                  \int_{r'\leq a}P_{l}(\mbox{\boldmath${\hat{k}}\cdot{\hat{r}}'$})
                        P_{n}(\mbox{\boldmath${\hat{r}}\cdot{\hat{r}}'$})
                             j_{l}(k_{1}r')j_{n}(kr_{<})h^{(1)}_{n}(kr_{>})d\mbox{\boldmath$r'$}\,.  
\end{eqnarray}
\newpage
\mbox{   }
\newline
\newline
\newline
It is necessary here to readjust the radial integration at $r'=r<a$ in conformity with the fact that the
decomposition (28) for the Green's function is sensitive to the relative order among magnitudes $r,$ $r'.$
After some rearrangement the right member of (32) becomes
\begin{eqnarray}
\lefteqn{
  \frac{\pi}{2i}\sqrt{\frac{k}{k_{1}}\,}\sum_{l=0}^{\infty}A^{<}_{l}P_{l}(\mbox{\boldmath${\hat{k}}\cdot{\hat{r}}$})
 \left[\rule{0mm}{6mm}r\sqrt{kk_{1}\,}j_{l}(k_{1}r)\left\{\rule{0mm}{5mm}
     J_{l+\frac{1}{2}}(kr)H^{(1)'}_{l+\frac{1}{2}}(kr)-J_{l+\frac{1}{2}}'(kr)H^{(1)}_{l+\frac{1}{2}}(kr)
      \right\}\rule{2cm}{0mm}\right.  }  \nonumber \\
 &    & \left.\rule{0mm}{6mm}\rule{3.8cm}{0mm}+aj_{l}(kr)\left\{\rule{0mm}{5mm}
   k_{1}J_{l+\frac{1}{2}}'(k_{1}a)H^{(1)}_{l+\frac{1}{2}}(ka)-kJ_{l+\frac{1}{2}}(k_{1}a)H^{(1)'}_{l+\frac{1}{2}}(ka)
     \right\}\right]\,.    
\end{eqnarray}
The coefficients of the various functions $j_{l}(k_{1}r)$ in expression (33) are seen to contain a well known
Wronskian with value $2i/\pi kr,$ which is independent of index $l.$  This value assures exact balance between
terms involving the functions $j_{l}(k_{1}r)$ in (32).  A requirement that the remaining series in
$P_{l}(\mbox{\boldmath${\hat{k}}\cdot{\hat{r}}$})j_{l}(kr)$ agree term-by-term supplies at length the fully explicit evaluations
\begin{eqnarray}
A^{<}_{l} & = & \frac{2i^{l-1}}{\pi a}\sqrt{\frac{k_{1}}{k}\,}(2l+1)\left(\rule{0mm}{6mm}
  k_{1}J_{l+\frac{1}{2}}'(k_{1}a)H^{(1)}_{l+\frac{1}{2}}(ka)-kJ_{l+\frac{1}{2}}(k_{1}a)H^{(1)'}_{l+\frac{1}{2}}(ka)\right)^{-1}\,.
\end{eqnarray}
It is then easily verified that (31) and (34) together point toward results in complete accord with those which follow from
the usual demand that both $\psi$ and its radial derivative be continuous at the potential barrier/well boundary.

    We encounter here a recurring pattern of great power, a decisive {\em{leitmotif}} in the inner/outer solution to inhomogeneous
equations such as (8).  As already embodied in expression (33) and Eq. (34), this pattern asserts itself within the scattering interior, wherein it provides an exact cancellation, left and right, of the total field (radial behavior guided by functions $j_{l}(k_{1}r)$).
Interior expansion coefficients $A^{<}_{l}$ are gotten then explicitly and without further ado from a balance along
{\em{exterior}} eigenfunctions (radial behavior guided by the $j_{l}(kr)$).  The exterior coefficients $A^{>}_{l}$
follow then directly from the $A^{<}_{l}$ as a sort of afterthought, as in (31), the link having been forged beforehand by
insisting that the integral equation be satisfied on the scattering exterior.

    This pattern of interior field cancellation recurs over and over again in applications considerably
more robust than the present quantum mechanical, strictly scalar setting.  Such applications are listed in a concluding
section.  We pause here only to stress the analytic economy of this two-tier, completely explicit solution process.   
\subsection{Solution based on the homogeneous form (11)}
     Scattering solution via the homogeneous form (11) follows similar lines, although the computations are consideraly
more laborious since now the necessity to segment radial integrations into distinct patterns above
and below observation
radius $r$ intrudes not only for $r<a$ but also on scatterer exterior $r>a.$   
Isolation of coefficients $A^{>}_{l}$ and $A^{<}_{l}$ proceeds as before by way of a term-by-term series identification,
as follows.
\newpage
\mbox{   }
\newline

     One begins by setting down as master relations corresponding respectively to (29) for $r>a$
\begin{eqnarray}
\lefteqn{
\sum_{l=0}^{\infty}\left(\rule{0mm}{5mm}A^{>}_{l}h_{l}^{(1)}(kr)+i^{l}(2l+1)j_{l}(kr)\right)
	           P_{l}(\mbox{\boldmath${\hat{k}}\cdot{\hat{r}}$})  = }  \nonumber  \\
&  &  \frac{mk_{0}V_{1}}{2\hbar^{2}\pi i}\sum_{l,n=0}^{\infty}A^{<}_{l}(2n+1)
h^{(1)}_{n}(k_{0}r)\int_{r'\leq a}P_{l}(\mbox{\boldmath${\hat{k}}\cdot{\hat{r}}'$})
P_{n}(\mbox{\boldmath${\hat{r}}\cdot{\hat{r}}'$})j_{l}(k_{1}r')j_{n}(k_{0}r')d\mbox{\boldmath$r'$} +   \\
&  &  \frac{mk_{0}V_{\infty}}{2\hbar^{2}\pi i}\sum_{l,n=0}^{\infty}(2n+1)
\int_{r'> a}P_{l}(\mbox{\boldmath${\hat{k}}\cdot{\hat{r}}'$})P_{n}(\mbox{\boldmath${\hat{r}}\cdot{\hat{r}}'$})\times \nonumber   \\
&  & \rule{3.3cm}{0mm} \times \left(\rule{0mm}{5mm}A^{>}_{l}h_{l}^{(1)}(kr')+i^{l}(2l+1)j_{l}(kr')\right)
j_{n}(k_{0}r_{<})h^{(1)}_{n}(k_{0}r_{>})d\mbox{\boldmath$r'$} \,,  \nonumber
\end{eqnarray}
and to (32) for $0\leq r \leq a$
\begin{eqnarray}
\lefteqn{
	\sum_{l=0}^{\infty}A^{<}_{l}j_{l}(k_{1}r)P_{l}(\mbox{\boldmath${\hat{k}}\cdot{\hat{r}}$})  = }  \nonumber  \\
&  &  \frac{mk_{0}V_{1}}{2\hbar^{2}\pi i}\sum_{l,n=0}^{\infty}A^{<}_{l}(2n+1)
\int_{r'\leq a}P_{l}(\mbox{\boldmath${\hat{k}}\cdot{\hat{r}}'$})
P_{n}(\mbox{\boldmath${\hat{r}}\cdot{\hat{r}}'$})j_{l}(k_{1}r')j_{n}(k_{0}r_{<})h^{(1)}_{n}(k_{0}r_{>})d\mbox{\boldmath$r'$}+   \\
&  &  \frac{mk_{0}V_{\infty}}{2\hbar^{2}\pi i}\sum_{l,n=0}^{\infty}(2n+1)j_{n}(k_{0}r)
\int_{r'> a}P_{l}(\mbox{\boldmath${\hat{k}}\cdot{\hat{r}}'$})P_{n}(\mbox{\boldmath${\hat{r}}\cdot{\hat{r}}'$})\times  \nonumber \\
&  &  \rule{4.4cm}{0mm}\times\left(\rule{0mm}{5mm}A^{>}_{l}h_{l}^{(1)}(kr')+i^{l}(2l+1)j_{l}(kr')\right)
h^{(1)}_{n}(k_{0}r')d\mbox{\boldmath$r'$} \,.  \nonumber
\end{eqnarray}
Here we utilize of course the Green's function analog of (28) gotten under the replacement of wave number $k$ by $k_{0}.$
Since Eqs. (35)-(36) share a comparable level of complexity, it suffices perhaps to follow the reduction of just one of
them, say for definiteness (35).
\subsubsection{Analytic reductions}
     And so as before in connection with (29) an integration over the solid angle of $\mbox{\boldmath$r'$}$ in (35) diagonalizes
indices $l$ and $n$ with
\begin{eqnarray}
\lefteqn{
	\sum_{l=0}^{\infty}\left(\rule{0mm}{5mm}A^{>}_{l}h_{l}^{(1)}(kr)+i^{l}(2l+1)j_{l}(kr)\right)
	P_{l}(\mbox{\boldmath${\hat{k}}\cdot{\hat{r}}$})  = }  \nonumber  \\
&  &  \frac{2mk_{0}V_{1}}{\hbar^{2} i}\sum_{l=0}^{\infty}A^{<}_{l}	P_{l}(\mbox{\boldmath${\hat{k}}\cdot{\hat{r}}$})
h^{(1)}_{l}(k_{0}r)\int_{r'\leq a}j_{l}(k_{1}r')j_{l}(k_{0}r')r'^{\,2}dr' +    \\
&  &  \frac{2mk_{0}V_{\infty}}{\hbar^{2} i}\sum_{l=0}^{\infty}P_{l}(\mbox{\boldmath${\hat{k}}\cdot{\hat{r}}$})
\int_{r'> a} \left(\rule{0mm}{5mm}A^{>}_{l}h_{l}^{(1)}(kr')+i^{l}(2l+1)j_{l}(kr')\right)
j_{l}(k_{0}r_{<})h^{(1)}_{l}(k_{0}r_{>})r'^{\,2}dr'   \nonumber
\end{eqnarray}     
as an intermediate result.  Here the self-field is conveyed by coefficients $A_{l}^{>},$ and one can initiate its ultimate
left/right cancellation, akin to that revealed in (33), by writing, for the relevant portion of the right-hand side,
\newpage
\mbox{   }
\newline
\begin{eqnarray}
\lefteqn{
\frac{2mk_{0}V_{\infty}}{\hbar^{2} i}\sum_{l=0}^{\infty}A^{>}_{l}P_{l}(\mbox{\boldmath${\hat{k}}\cdot{\hat{r}}$})
\int_{r'> a}h_{l}^{(1)}(kr')j_{l}(k_{0}r_{<})h^{(1)}_{l}(k_{0}r_{>})r'^{\,2}dr' = } \nonumber \\
   &   & \rule{2cm}{0mm}  \frac{2mk_{0}V_{\infty}}{\hbar^{2} i}\sum_{l=0}^{\infty}A^{>}_{l}P_{l}(\mbox{\boldmath${\hat{k}}\cdot{\hat{r}}$})
   h^{(1)}_{l}(k_{0}r)\int_{a}^{r}h_{l}^{(1)}(kr')j_{l}(k_{0}r')r'^{\,2}dr' +          \\
   &   & \rule{2cm}{0mm}  \frac{2mk_{0}V_{\infty}}{\hbar^{2} i}\sum_{l=0}^{\infty}A^{>}_{l}P_{l}(\mbox{\boldmath${\hat{k}}\cdot{\hat{r}}$})
   j_{l}(k_{0}r)\int_{r}^{\infty}h_{l}^{(1)}(kr')h^{(1)}_{l}(k_{0}r')r'^{\,2}dr'\,.  \nonumber
\end{eqnarray}
Application of the Bessel/Lommel quadrature, previously cited, converts this into
{\small{
\begin{eqnarray}
\lefteqn{
	\frac{2mk_{0}V_{\infty}}{\hbar^{2} i}\sum_{l=0}^{\infty}A^{>}_{l}P_{l}(\mbox{\boldmath${\hat{k}}\cdot{\hat{r}}$})
	\int_{r'> a}h_{l}^{(1)}(kr')j_{l}(k_{0}r_{<})h^{(1)}_{l}(k_{0}r_{>})r'^{\,2}dr' = } \nonumber \\
&   & \frac{\pi mk_{0}V_{\infty}}{i\hbar^{2}(k_{0}^{2}-k^{2})\sqrt{kk_{0}\,}}
    \sum_{l=0}^{\infty}A^{>}_{l}P_{l}(\mbox{\boldmath${\hat{k}}\cdot{\hat{r}}$})h^{(1)}_{l}(k_{0}r)
\left[\rule{0mm}{7mm}\left.\rule{0mm}{6mm}r'\left(\rule{0mm}{5mm}kJ_{l+\frac{1}{2}}(k_{0}r')H^{(1)'}_{l+\frac{1}{2}}(kr')-
                k_{0}J'_{l+\frac{1}{2}}(k_{0}r')H^{(1)}_{l+\frac{1}{2}}(kr')\right)\right|_{r'=a}^{r'=r}\right]   \\
&   & +\frac{\pi mk_{0}V_{\infty}}{i\hbar^{2}(k_{0}^{2}-k^{2})\sqrt{kk_{0}\,}}
\sum_{l=0}^{\infty}A^{>}_{l}P_{l}(\mbox{\boldmath${\hat{k}}\cdot{\hat{r}}$})j_{l}(k_{0}r)
\left[\rule{0mm}{7mm}\left.\rule{0mm}{6mm}r'\left(\rule{0mm}{5mm}kH^{(1)}_{l+\frac{1}{2}}(k_{0}r')H^{(1)'}_{l+\frac{1}{2}}(kr')-
         k_{0}H^{(1)'}_{l+\frac{1}{2}}(k_{0}r')H^{(1)}_{l+\frac{1}{2}}(kr')\right)\right|_{r'=r}^{r'=\infty}\right]\,. \nonumber
\end{eqnarray} }  }\vspace{-6mm}

\parindent=0in
On taking note from (3) and (10) of the fact, akin to (30), that $k_{0}^{2}-k^{2}=2mV_{\infty}/\hbar^{2}$ we arrive at the next intermediate
stage
	{\small{
			\begin{eqnarray}
			\lefteqn{
				\frac{2mk_{0}V_{\infty}}{\hbar^{2} i}\sum_{l=0}^{\infty}A^{>}_{l}P_{l}(\mbox{\boldmath${\hat{k}}\cdot{\hat{r}}$})
				\int_{r'> a}h_{l}^{(1)}(kr')j_{l}(k_{0}r_{<})h^{(1)}_{l}(k_{0}r_{>})r'^{\,2}dr' = } \nonumber \\
			&   & \rule{2mm}{0mm}\frac{\pi^{3/2}}{2i\sqrt{2kr\,}}
			\sum_{l=0}^{\infty}A^{>}_{l}P_{l}(\mbox{\boldmath${\hat{k}}\cdot{\hat{r}}$})H^{(1)}_{l+\frac{1}{2}}(k_{0}r)
	\left[\rule{0mm}{7mm}\left.\rule{0mm}{6mm}r'\left(\rule{0mm}{5mm}kJ_{l+\frac{1}{2}}(k_{0}r')H^{(1)'}_{l+\frac{1}{2}}(kr')
	         - k_{0}J'_{l+\frac{1}{2}}(k_{0}r')H^{(1)}_{l+\frac{1}{2}}(kr')\right)\right|_{r'=a}^{r'=r}\right]   \\
			&   & +\frac{\pi^{3/2}}{2i\sqrt{2kr\,}}
			\sum_{l=0}^{\infty}A^{>}_{l}P_{l}(\mbox{\boldmath${\hat{k}}\cdot{\hat{r}}$})J_{l+\frac{1}{2}}(k_{0}r)
	\left[\rule{0mm}{7mm}\left.\rule{0mm}{6mm}r'\left(\rule{0mm}{5mm}kH^{(1)}_{l+\frac{1}{2}}(k_{0}r')H^{(1)'}_{l+\frac{1}{2}}(kr')
	        - k_{0}H^{(1)'}_{l+\frac{1}{2}}(k_{0}r')H^{(1)}_{l+\frac{1}{2}}(kr')\right)\right|_{r'=r}^{r'=\infty}\right]\,, \nonumber
			\end{eqnarray} }  }\vspace{-6mm}
		
\parindent=0in	
a stage which admits a more perspicuous rearrangement of terms so as to bring into focus once again a recurring appearance
by our most useful Bessel/Hankel Wronskian.  Thus
		\begin{eqnarray}
		\lefteqn{
			\frac{2mk_{0}V_{\infty}}{\hbar^{2} i}\sum_{l=0}^{\infty}A^{>}_{l}P_{l}(\mbox{\boldmath${\hat{k}}\cdot{\hat{r}}$})
			\int_{r'> a}h_{l}^{(1)}(kr')j_{l}(k_{0}r_{<})h^{(1)}_{l}(k_{0}r_{>})r'^{\,2}dr' = } \nonumber \\
		&   & \rule{2mm}{0mm}\frac{\pi^{3/2}}{2i\sqrt{2kr\,}}
		\sum_{l=0}^{\infty}A^{>}_{l}P_{l}(\mbox{\boldmath${\hat{k}}\cdot{\hat{r}}$})H^{(1)}_{l+\frac{1}{2}}(kr)k_{0}r
		\left\{\rule{0mm}{5mm}J_{l+\frac{1}{2}}(k_{0}r)H^{(1)'}_{l+\frac{1}{2}}(k_{0}r)-
		             J'_{l+\frac{1}{2}}(k_{0}r)H^{(1)}_{l+\frac{1}{2}}(k_{0}r)\right\}  \nonumber \\
		&   & -\frac{\pi^{3/2}}{2i\sqrt{2kr\,}}
		\sum_{l=0}^{\infty}A^{>}_{l}P_{l}(\mbox{\boldmath${\hat{k}}\cdot{\hat{r}}$})H^{(1)}_{l+\frac{1}{2}}(k_{0}r)a
		\left\{\rule{0mm}{5mm}kJ_{l+\frac{1}{2}}(k_{0}a)H^{(1)'}_{l+\frac{1}{2}}(ka)-
		     k_{0}J_{l+\frac{1}{2}}'(k_{0}a)H^{(1)}_{l+\frac{1}{2}}(ka)\right\}      \\
		& = & \frac{i\pi a}{2}\sqrt{\frac{k_{0}}{k}\,}
		\sum_{l=0}^{\infty}A^{>}_{l}P_{l}(\mbox{\boldmath${\hat{k}}\cdot{\hat{r}}$})h^{(1)}_{l}(k_{0}r)
		\left\{\rule{0mm}{5mm}kJ_{l+\frac{1}{2}}(k_{0}a)H^{(1)'}_{l+\frac{1}{2}}(ka)-
		k_{0}J_{l+\frac{1}{2}}'(k_{0}a)H^{(1)}_{l+\frac{1}{2}}(ka)\right\}  \nonumber \\
		&   &   \rule{1.1cm}{0mm} +
		\sum_{l=0}^{\infty}A^{>}_{l}P_{l}(\mbox{\boldmath${\hat{k}}\cdot{\hat{r}}$})h^{(1)}_{l}(kr) \nonumber 
		\end{eqnarray}
\newpage
\mbox{   }
\newline
\newline
\newline
which leads of course to an exact cancellation, {\em{vis-\`{a}-vis}} Eq. (37), of the self-field
$\sum_{l=0}^{\infty}A^{>}_{l}P_{l}(\mbox{\boldmath${\hat{k}}\cdot{\hat{r}}$})h^{(1)}_{l}(kr).$
\parindent=0.5in

    The remaining ingredients of Eq. (37) are processed in similar fashion.  The first term on the right is the
easiest.  It gives
\begin{eqnarray}
\lefteqn{
\frac{2mk_{0}V_{1}}{\hbar^{2} i}\sum_{l=0}^{\infty}A^{<}_{l}P_{l}(\mbox{\boldmath${\hat{k}}\cdot{\hat{r}}$})
h^{(1)}_{l}(k_{0}r)\int_{r'\leq a}j_{l}(k_{1}r')j_{l}(k_{0}r')r'^{\,2}dr' = }  \nonumber \\
  &  & \rule{1mm}{0mm} \frac{\pi m a V_{1}}{\hbar^{2} i (k_{0}^{2}-k_{1}^{2})}\sqrt{\frac{k_{0}}{k_{1}}\,}
         \sum_{l=0}^{\infty}A^{<}_{l}P_{l}(\mbox{\boldmath${\hat{k}}\cdot{\hat{r}}$})h^{(1)}_{l}(k_{0}r)
  \left\{\rule{0mm}{5mm} k_{1}J_{l+\frac{1}{2}}(k_{0}a)J_{l+\frac{1}{2}}'(k_{1}a)-
                         k_{0}J_{l+\frac{1}{2}}'(k_{0}a)J_{l+\frac{1}{2}}(k_{1}a)\right\} \nonumber   \\
  & = & -\frac{i\pi a}{2}\sqrt{\frac{k_{0}}{k_{1}}\,}
    \sum_{l=0}^{\infty}A^{<}_{l}P_{l}(\mbox{\boldmath${\hat{k}}\cdot{\hat{r}}$})h^{(1)}_{l}(k_{0}r)
  \left\{\rule{0mm}{5mm} k_{1}J_{l+\frac{1}{2}}(k_{0}a)J_{l+\frac{1}{2}}'(k_{1}a)-
            k_{0}J_{l+\frac{1}{2}}'(k_{0}a)J_{l+\frac{1}{2}}(k_{1}a)\right\}\,,                           	                   
\end{eqnarray}
wherein $V_{1}$ as an overt factor vanishes by virtue of Eqs. (10) and (26).

      For the final term on the right in (37), involving the incoming excitation, we are forced to emulate steps
(38)-(41).  Thus
\begin{eqnarray}
\lefteqn{
	\frac{2mk_{0}V_{\infty}}{\hbar^{2} i}\sum_{l=0}^{\infty}i^{l}(2l+1)P_{l}(\mbox{\boldmath${\hat{k}}\cdot{\hat{r}}$})
	\int_{r'> a}j_{l}(kr')j_{l}(k_{0}r_{<})h^{(1)}_{l}(k_{0}r_{>})r'^{\,2}dr' = } \nonumber \\
&   & \rule{2cm}{0mm}  \frac{2mk_{0}V_{\infty}}{\hbar^{2} i}\sum_{l=0}^{\infty}
   i^{l}(2l+1)P_{l}(\mbox{\boldmath${\hat{k}}\cdot{\hat{r}}$})
h^{(1)}_{l}(k_{0}r)\int_{a}^{r}j_{l}(kr')j_{l}(k_{0}r')r'^{\,2}dr' +          \\
&   & \rule{2cm}{0mm}  \frac{2mk_{0}V_{\infty}}{\hbar^{2} i}\sum_{l=0}^{\infty}
i^{l}(2l+1)P_{l}(\mbox{\boldmath${\hat{k}}\cdot{\hat{r}}$})
j_{l}(k_{0}r)\int_{r}^{\infty}j_{l}(kr')h^{(1)}_{l}(k_{0}r')r'^{\,2}dr'\,.  \nonumber
\end{eqnarray}
Without further ado there follows next a cascade of three steps.  First
{\small{
		\begin{eqnarray}
		\lefteqn{
			\frac{2mk_{0}V_{\infty}}{\hbar^{2} i}\sum_{l=0}^{\infty}
			i^{l}(2l+1)P_{l}(\mbox{\boldmath${\hat{k}}\cdot{\hat{r}}$})
			\int_{r'> a}j_{l}(kr')j_{l}(k_{0}r_{<})h^{(1)}_{l}(k_{0}r_{>})r'^{\,2}dr' = }  \\
		&   & \frac{\pi mk_{0}V_{\infty}}{i\hbar^{2}(k_{0}^{2}-k^{2})\sqrt{kk_{0}\,}}
		\sum_{l=0}^{\infty}i^{l}(2l+1)P_{l}(\mbox{\boldmath${\hat{k}}\cdot{\hat{r}}$})
		 h^{(1)}_{l}(k_{0}r)
		\left[\rule{0mm}{7mm}\left.\rule{0mm}{6mm}r'\left(\rule{0mm}{5mm}kJ_{l+\frac{1}{2}}(k_{0}r')J'_{l+\frac{1}{2}}(kr')-
		k_{0}J'_{l+\frac{1}{2}}(k_{0}r')J_{l+\frac{1}{2}}(kr')\right)\right|_{r'=a}^{r'=r}\right]  \nonumber \\
		&   & +\frac{\pi mk_{0}V_{\infty}}{i\hbar^{2}(k_{0}^{2}-k^{2})\sqrt{kk_{0}\,}}
		\sum_{l=0}^{\infty}i^{l}(2l+1)P_{l}(\mbox{\boldmath${\hat{k}}\cdot{\hat{r}}$})j_{l}(k_{0}r)
		\left[\rule{0mm}{7mm}\left.\rule{0mm}{6mm}r'\left(\rule{0mm}{5mm}kH^{(1)}_{l+\frac{1}{2}}(k_{0}r')J'_{l+\frac{1}{2}}(kr')-
		k_{0}H^{(1)'}_{l+\frac{1}{2}}(k_{0}r')J_{l+\frac{1}{2}}(kr')\right)\right|_{r'=r}^{r'=\infty}\right]\,, \nonumber
		\end{eqnarray} }  }
\vspace{-9mm}

\parindent=0in
then
\vspace{2mm}
{\small{
		\begin{eqnarray}
		\lefteqn{
			\frac{2mk_{0}V_{\infty}}{\hbar^{2} i}\sum_{l=0}^{\infty}i^{l}(2l+1)P_{l}(\mbox{\boldmath${\hat{k}}\cdot{\hat{r}}$})
			\int_{r'> a}j_{l}(kr')j_{l}(k_{0}r_{<})h^{(1)}_{l}(k_{0}r_{>})r'^{\,2}dr' = } \nonumber \\
		&   & \rule{2mm}{0mm}\frac{\pi^{3/2}}{2i\sqrt{2kr\,}}
		\sum_{l=0}^{\infty}i^{l}(2l+1)P_{l}(\mbox{\boldmath${\hat{k}}\cdot{\hat{r}}$})H^{(1)}_{l+\frac{1}{2}}(k_{0}r)
		\left[\rule{0mm}{7mm}\left.\rule{0mm}{6mm}r'\left(\rule{0mm}{5mm}kJ_{l+\frac{1}{2}}(k_{0}r')J'_{l+\frac{1}{2}}(kr')
		- k_{0}J'_{l+\frac{1}{2}}(k_{0}r')J_{l+\frac{1}{2}}(kr')\right)\right|_{r'=a}^{r'=r}\right]   \\
		&   & +\frac{\pi^{3/2}}{2i\sqrt{2kr\,}}
		\sum_{l=0}^{\infty}i^{l}(2l+1)P_{l}(\mbox{\boldmath${\hat{k}}\cdot{\hat{r}}$})J_{l+\frac{1}{2}}(k_{0}r)
		\left[\rule{0mm}{7mm}\left.\rule{0mm}{6mm}r'\left(\rule{0mm}{5mm}kH^{(1)}_{l+\frac{1}{2}}(k_{0}r')J'_{l+\frac{1}{2}}(kr')
		- k_{0}H^{(1)'}_{l+\frac{1}{2}}(k_{0}r')J_{l+\frac{1}{2}}(kr')\right)\right|_{r'=r}^{r'=\infty}\right]\,, \nonumber
		\end{eqnarray} }  }
\parindent=0in
\newpage
\mbox{   }
\newline
\newline
and finally
		\begin{eqnarray}
		\lefteqn{
			\frac{2mk_{0}V_{\infty}}{\hbar^{2} i}\sum_{l=0}^{\infty}i^{l}(2l+1)P_{l}(\mbox{\boldmath${\hat{k}}\cdot{\hat{r}}$})
			\int_{r'> a}j_{l}(kr')j_{l}(k_{0}r_{<})h^{(1)}_{l}(k_{0}r_{>})r'^{\,2}dr' = } \nonumber \\
		&   & \rule{2mm}{0mm}\frac{\pi^{3/2}}{2i\sqrt{2kr\,}}
		\sum_{l=0}^{\infty}i^{l}(2l+1)P_{l}(\mbox{\boldmath${\hat{k}}\cdot{\hat{r}}$})J_{l+\frac{1}{2}}(kr)k_{0}r
		\left\{\rule{0mm}{5mm}J_{l+\frac{1}{2}}(k_{0}r)H^{(1)'}_{l+\frac{1}{2}}(k_{0}r)-
		J'_{l+\frac{1}{2}}(k_{0}r)H^{(1)}_{l+\frac{1}{2}}(k_{0}r)\right\}  \nonumber \\
		&   & -\frac{\pi^{3/2}}{2i\sqrt{2kr\,}}
		\sum_{l=0}^{\infty}i^{l}(2l+1)P_{l}(\mbox{\boldmath${\hat{k}}\cdot{\hat{r}}$})H^{(1)}_{l+\frac{1}{2}}(k_{0}r)a
		\left\{\rule{0mm}{5mm}kJ_{l+\frac{1}{2}}(k_{0}a)J'_{l+\frac{1}{2}}(ka)-
		k_{0}J_{l+\frac{1}{2}}'(k_{0}a)J_{l+\frac{1}{2}}(ka)\right\}      \\
		& = & \frac{i\pi a}{2}\sqrt{\frac{k_{0}}{k}\,}
		\sum_{l=0}^{\infty}i^{l}(2l+1)P_{l}(\mbox{\boldmath${\hat{k}}\cdot{\hat{r}}$})h^{(1)}_{l}(k_{0}r)
		\left\{\rule{0mm}{5mm}kJ_{l+\frac{1}{2}}(k_{0}a)J'_{l+\frac{1}{2}}(ka)-
		k_{0}J_{l+\frac{1}{2}}'(k_{0}a)J_{l+\frac{1}{2}}(ka)\right\} \,, \nonumber \\
		&   &   \rule{1.1cm}{0mm} +
		\sum_{l=0}^{\infty}i^{l}(2l+1)P_{l}(\mbox{\boldmath${\hat{k}}\cdot{\hat{r}}$})j_{l}(kr)\,, \nonumber 
		\end{eqnarray} 
which assures yet another self-cancellation against the left side in (37), but accompanied now by a drive contribution
itself proportional to $h^{(1)}_{l}(k_{0}r),$ just as in (41) preceding.\footnote{In arriving at (46) we have glibly
	discarded the Lommel quadrature contribution at $r'=\infty$ in (45).  Such omission is predicated on having wave number
	$k_{0}$ in Green's function (28) endowed with a small positive imaginary component, already announced beneath
	Footnote 3, admittedly a {\em{deus ex machina}} artifice which, in physical terms, implies some sort of particle
	disappearance, evaporation, if you would, on its outward flight to infinity, a concept altogether foreign to Schr{\"{o}}dinger's
	equation as normally received.  In the absence of such an imaginary component the Lommel contribution at infinity is
	{\em{``infinite oscillatory,''}} and thus mathematically indigestible.  It is perhaps better to simply accept, and learn
	to live with the outbound {\em{``particle probability dissipation.''}}  Of course, in view of the still earlier Footnote 1,
	there does hover the redemptive possibility of $k_{0}$ actually being, in its own right, positive imaginary at any magnitude.
    Similar anxieties regarding an undeclared quadrature omission at infinity when passing from (40) to (41) are softened even
    more decisively by the fact that one confronts there a pair of outgoing Hankel functions, each one of which is protected by a dissipative mechanism.}	
\parindent=0.5in

    Lastly, on putting together the remaining pieces from (41), (42), and (46) we arrive at
\begin{eqnarray}
\lefteqn{
\sum_{l=0}^{\infty}P_{l}(\mbox{\boldmath${\hat{k}}\cdot{\hat{r}}$})h^{(1)}_{l}(k_{0}r)\left[\rule{0mm}{7mm}\rule{1mm}{0mm}
A^{>}_{l}\sqrt{\frac{1}{k}\,}\left\{\rule{0mm}{5mm}kJ_{l+\frac{1}{2}}(k_{0}a)H^{(1)'}_{l+\frac{1}{2}}(ka)-
k_{0}J_{l+\frac{1}{2}}'(k_{0}a)H^{(1)}_{l+\frac{1}{2}}(ka)\right\}-\right. } \nonumber \\
  &    & \rule{2.7cm}{0mm} \left.\rule{0mm}{7mm}
     A^{<}_{l}\sqrt{\frac{1}{k_{1}}\,}     
     \left\{\rule{0mm}{5mm} k_{1}J_{l+\frac{1}{2}}(k_{0}a)J_{l+\frac{1}{2}}'(k_{1}a)-
     k_{0}J_{l+\frac{1}{2}}'(k_{0}a)J_{l+\frac{1}{2}}(k_{1}a)\right\}\right. +         \\
  &    & \left.\rule{0mm}{7mm}\rule{2.0cm}{0mm}
   i^{l}(2l+1)\sqrt{\frac{1}{k}\,}
   \left\{\rule{0mm}{5mm}kJ_{l+\frac{1}{2}}(k_{0}a)J'_{l+\frac{1}{2}}(ka)-
   k_{0}J_{l+\frac{1}{2}}'(k_{0}a)J_{l+\frac{1}{2}}(ka)\right\}\rule{1mm}{0mm} \right]\,=\,0\,. \nonumber
\end{eqnarray} 
a global statement which splinters of course into an infinity of subsidiary demands			     	                 
		\begin{eqnarray}
		\lefteqn{			
			A^{>}_{l}\sqrt{\frac{1}{k}\,}\left\{\rule{0mm}{5mm}kJ_{l+\frac{1}{2}}(k_{0}a)H^{(1)'}_{l+\frac{1}{2}}(ka)-
			k_{0}J_{l+\frac{1}{2}}'(k_{0}a)H^{(1)}_{l+\frac{1}{2}}(ka)\right\}- } \nonumber \\
		&    & \rule{1.5cm}{0mm} 
		A^{<}_{l}\sqrt{\frac{1}{k_{1}}\,}     
		\left\{\rule{0mm}{5mm} k_{1}J_{l+\frac{1}{2}}(k_{0}a)J_{l+\frac{1}{2}}'(k_{1}a)-
		k_{0}J_{l+\frac{1}{2}}'(k_{0}a)J_{l+\frac{1}{2}}(k_{1}a)\right\}  +         \\
		&    &   \rule{2.5cm}{0mm}
		i^{l}(2l+1)\sqrt{\frac{1}{k}\,}
		\left\{\rule{0mm}{5mm}kJ_{l+\frac{1}{2}}(k_{0}a)J'_{l+\frac{1}{2}}(ka)-
		k_{0}J_{l+\frac{1}{2}}'(k_{0}a)J_{l+\frac{1}{2}}(ka)\right\}  \,=\,0\,. \nonumber
		\end{eqnarray} 
holding good $\forall$ $l\geq 0.$
\newpage
\mbox{   }
\newline

     A reduction of comparable complexity is likewise available for Eq. (36), its quadrature disjunction at
$r'=r<a$ affecting only terms from the first summation on its right.  It permits us to adjoin to (48) the independent
statements
\begin{eqnarray}
\lefteqn{			
	A^{>}_{l}\sqrt{\frac{1}{k}\,}\left\{\rule{0mm}{5mm}kH^{(1)}_{l+\frac{1}{2}}(k_{0}a)H^{(1)'}_{l+\frac{1}{2}}(ka)-
	k_{0}H^{(1)'}_{l+\frac{1}{2}}(k_{0}a)H^{(1)}_{l+\frac{1}{2}}(ka)\right\}- } \nonumber \\
&    & \rule{1.5cm}{0mm} 
A^{<}_{l}\sqrt{\frac{1}{k_{1}}\,}     
\left\{\rule{0mm}{5mm} k_{1}H^{(1)}_{l+\frac{1}{2}}(k_{0}a)J_{l+\frac{1}{2}}'(k_{1}a)-
k_{0}H^{(1)'}_{l+\frac{1}{2}}(k_{0}a)J_{l+\frac{1}{2}}(k_{1}a)\right\}  +         \\
&    &   \rule{2.5cm}{0mm}
i^{l}(2l+1)\sqrt{\frac{1}{k}\,}
\left\{\rule{0mm}{5mm}kH^{(1)}_{l+\frac{1}{2}}(k_{0}a)J'_{l+\frac{1}{2}}(ka)-
k_{0}H^{(1)'}_{l+\frac{1}{2}}(k_{0}a)J_{l+\frac{1}{2}}(ka)\right\}  \,=\,0 \nonumber
\end{eqnarray}
holding good as before $\forall$ $l\geq 0.$
\vspace{-2mm}
\subsubsection{Final analytic synopsis}
\vspace{-0mm}
     Equations (48)-(49) are more concisely rendered in a matricial format governing coefficient vectors
$A_{l}$ defined by
\begin{equation}
A_{l}  =  \left[\begin{array}{c}
                      A^{>}_{l} \\
                                \\
                      A^{<}_{l}
                \end{array} \right] \,.
\end{equation}
One then finds for all indices $l\geq 0$ that
\begin{eqnarray}
M_{l}\left(\rule{0mm}{4.0mm}N_{l}A_{l}-B_{l}\right) & = & 0\,,
\end{eqnarray}                
with matrices
\begin{equation}
M_{l} = \left[\rule{0mm}{1.1cm}\begin{array}{lr}
                       J_{l+\frac{1}{2}}(k_{0}a)         &   -k_{0}J_{l+\frac{1}{2}}'(k_{0}a)     \\
                                                                                                  \\
                       H_{l+\frac{1}{2}}^{(1)}(k_{0}a)   &   -k_{0}H_{l+\frac{1}{2}}^{(1)'}(k_{0}a) 
                    \end{array}  \right]\,,
\end{equation}
\newline                     
\begin{equation}
N_{l} = \left[\rule{0mm}{1.1cm}\begin{array}{ll}
                      -\sqrt{k}H_{l+\frac{1}{2}}^{(1)'}(ka)    &     \sqrt{k_{1}}J_{l+\frac{1}{2}}'(k_{1}a)     \\
                                                                                                                \\
                      -H_{l+\frac{1}{2}}^{(1)}(ka)/\sqrt{k}    &     J_{l+\frac{1}{2}}(k_{1}a)/\sqrt{k_{1}} 
\end{array}  \right]\,,
\end{equation}                                                                   
and
\begin{equation}
B_{l} = \frac{i^{l}(2l+1)}{\sqrt{k\,}}\left[\begin{array}{c}
	                          kJ_{l+\frac{1}{2}}^{'}(ka)  \\
	                                                      \\
	                          J_{l+\frac{1}{2}}(ka)
	                          \end{array} \right]\,.                            
\end{equation}
Our ubiquitous Wronskian provides the nonvanishing evaluation $\det M_{l} = 2/\pi ia \neq 0 \,,$
and hence assures that the matrix factor $M_{l}$ in (51) is in fact irrelevant and may be freely detached.
The ensuing requirement that the remaining structure $N_{l}A_{l} - B_{l}$ on the left in (51)
vanish then easily leads to results consonant in every detail with those contained within (31) and (34).
The stated Wronskian participates yet again when agreement is first sought on behalf of $A^{<}_{l}.$
\newpage
\mbox{   }
\newline

    Indeed, the two-by-two linear system $N_{l}A_{l}-B_{l}=0$ is found to be identical,\footnote{Save for an
inessential row interchange and an overall minus sign.} as the revealed harmony
of $A_{l}$ outcomes now clearly requires it to be, with
that which expresses continuity at barrier/well interface $r=a$ for both $\psi$ and its radial derivative
$\partial\psi/\partial r.$     

     One draws much reassurance from this agreement, since it provides clear, specific evidence of an otherwise globally
proclaimed duality.  On the other hand, since so much turgid algebra is stirred up along this path, one would
hesitate to advocate on behalf of the homogeneous version as an efficient computational tool.
\vspace{-8mm}
\section{An {\textit{aperc\c{u}}} on h/i integral equations in scattering theory}
\vspace{-4mm}
Sporadic sightings of integral equations, both homogeneous and inhomogeneous, have been encountered in the scattering
literature, both electromagnetic and quantum mechanical.  We can leave to the side the entire subgenre of Born
and allied, Neumann-like iterative approximation schemes which figure so prominently in quantum field theory
and are abundantly documented.

     An early example of a homogeneous integral equation in the service of electromagnetics can be found in [{\bf{1}}],
offered therein as a vehicle for field/eigenfrequency determination in a cavity which encloses a dielectric object.
Such offer, alas, is bereft of any further examples as to its use.

    On an entirely different tack, now in solid state physics, there appeared the exquisitely beautiful paper of
Saxon and Hutner [{\bf{2}}], which extended and superseded the foundational Kronig-Penney treatment of
one-dimensional quantum mechanical electron movement in a periodic crystal lattice.  The Schr{\"{o}}dinger
equation was recast there under a homogeneous integral equation form, the lattice obstacles being modeled via
Dirac deltas and with the Bloch {\em{ansatz}} firmly in mind.

    The Saxon-Hutner one-dimensional achievement was soon to be eclipsed by the full, three-dimension-
al method
attributed to Kohn, Korringa, and Rostoker, the so-called KKR method [{\bf{3}}], so named after its originators.
Their primary contributions are found in [{\bf{4-5}}].     

    The homogeneous/inhomogeneous duality now on view seems first to have surfaced quite some time ago in a short    
note [{\bf{6}}] devoted to electromagnetic scattering by dielectric obstacles.  It turned out then that the
vector attribute of the electromagnetic field was an inessential complication and could easily be accommodated on
the road to duality demonstration.  On the other hand, the brevity of that note, and the paucity of detailed
development at that time, precluded the setting out of any specific application, such as that of the present
spherical scatterer.

     Over the ensuing decades occasional use was made to great effect of the self-field cancellation
phenomenon as a bridge to concrete solutions.  A first indication thereof appeared in [{\bf{7}}]
and dealt with electromagnetic reflection from and transmission into a lossy dielectric
half space.\footnote{It should perhaps be mentioned in passing that the self-field cancellation phenomenon is
reminiscent of the so called Ewald-Oseen extinction viewpoint encountered in optics.  We, however, refrain from
belaboring its physically obvious association with the concept of an incoming {\em{versus}} an internal field
exchange, and view it simply as a useful, systematic route to problem solution.}
\newpage
\mbox{   }
\newline

     Self-field cancellation showed its great power even in a magnetized plasma environment wherein dielectric
properties had risen to a tensorial level.  A pattern of such applications, originating in work by the undersigned,
can be found in [{\bf{8-9}}].

     Several decades later, aided now by Laplace transformation, a still broader arena of time-dependent scattering
was found to fit easily beneath a self-field cancellation purview [{\bf{10}}].  And lastly, more recent self-field 
cancellation efforts have yielded a fully vectorial solution for reflection from and penetration across a dielectric
slab [{\bf{11}}], and an inroad at least into the heretofore intractable problem of electromagnetic diffraction by
a dielectric wedge [{\bf{12}}].
\section{References}
\parindent=0in 

[1]	Julian Schwinger (1943).  {\bf{The Theory of Obstacles in Resonant Cavities and Waveguides}}, MIT Radiation Laboratory
Report 43-34, p. 16, eqs. (58)-(59).

[2]	D. S. Saxon, R. A. Hutner (1949).  {\bf{Some Electronic Properties of a One-dimensional Crystal Model}}, Philips
Research Reports, No. 4, pp. 82-122.

[3]	J. M. Ziman (1972).  {\bf{Principles of the Theory of Solids}}, {\em{Second Edition}}, Cambridge University Press,
pp. 106-108.

[4]	 J. Korringa (1947).  {\bf{On the calculation of the energy of a Bloch wave in a metal}}, Physica, Vol. XIII, Nos.
6–7, pp. 392-400.

[5]	W. Kohn, N. Rostoker (1954).  {\bf{Solution of the Schr{\"{o}}dinger Equation in Periodic Lattices with an Application to
Metallic Lithium}}, Phys. Rev., Vol. 94, No. 5, pp. 1111-1120.

[6]	J. Grzesik (1966).  {\bf{Note on Homogeneous and Inhomogeneous Integral Equations in the Theory of Electromagnetic
Scattering by Dielectric Obstacles}}, Proc. IEEE, Vol. 54, No. 12, pp. 2028-2029.

[7]	J. Grzesik (1980).  {\bf{Field Matching through Volume Suppression}}, IEE Proc., Vol. 127, Pt. H, No. 1, pp. 20-26.

[8]	Leif Ulstrup and Jan Grzesik (1984).  {\bf{Vector Green's Function Coil Code:  Numerical Results for a Half-Turn Loop
		in a Circular Waveguide}}, Twenty-sixth Annual Meeting, APS Division of Plasma Physics, Boston.

[9]	Hiroshi Agravante (1987).  {\bf{Plasma Separation Process:  Vector Green's Function Coil Code User's Manual}},
TRW Space and Technology Group, One Space Park, Redondo Beach, CA 90278, Report PSP-R1-1307.
\newpage
\mbox{   }
\newline
\newline
\newline
[10]	J. A. Grzesik (2007).  {\bf{EM Pulse Transit across a Uniform Dielectric Slab}}, $7^{{\rm{th}}}$ Workshop on
Computational Electromagnetics in Time-Domain, Perugia, Italy.

[11]	J. A. Grzesik (2018).  {\bf{Dielectric Slab Reflection/Transmission as a Self-Consistent Radiation Phenomenon}},     	     	 
Progress in Electromagnetics Research B, Vol. 82, pp. 31-48.

[12]	J. A. Grzesik (2019).  {\bf{Dielectric Wedge Scattering: An Analytic Inroad}}, Progress in Electromagnetics Research B,
Vol. 84, pp. 43-60.


\end{document}